\documentclass[twoside]{article}
\usepackage{amsmath}
\usepackage{amssymb}
\usepackage[mathscr]{eucal}
\usepackage{verbatim}
\usepackage{theorem}
\usepackage[hyphens]{url}

\newtheorem{theorem}{Theorem}[section]
\newtheorem{corollary}[theorem]{Corollary}
\newtheorem{proposition}[theorem]{Proposition}

{\theorembodyfont{\rmfamily}

}

{\theorembodyfont{\rmfamily}
    \newtheorem{example}[theorem]{Example}
    \newtheorem{remark}[theorem]{Remark}
}

\newenvironment{prooff}{\noindent{\em Proof.}}{$\blacksquare$\medskip}

\newcommand{\bsone}{\boldsymbol{1}}
\newcommand{\bstwo}{\boldsymbol{2}}
\newcommand{\bsl}{\boldsymbol{l}}

\newcommand{\mcJ}{\mathcal{J}}
\newcommand{\smcJ}{\text{\raisebox{0.15ex}{\small$\mathcal{J}$}}}
\newcommand{\mcK}{\mathcal{K}}
\newcommand{\mcL}{\mathcal{L}}

\newcommand{\mcV}{\mathcal{V}}
\newcommand{\mcW}{\mathcal{W}}

\newcommand{\mbR}{\mathbb{R}}
\newcommand{\mbS}{\mathbb{S}}
\newcommand{\mbZ}{\mathbb{Z}}

\newcommand{\Expec}{\mathcal{E}}
\newcommand{\Disper}{\mathcal{D}}

\newcommand{\eqval}{~\Leftrightarrow~}

\DeclareMathOperator{\Prob}{Pr}

\DeclareMathOperator{\rng}{rank}

\DeclareMathOperator{\Supp}{Supp}

\DeclareMathOperator{\Lin}{Lin}

\newcommand{\xaff}{\normalsize Department of Statistics,
                               Brigham Young University,
                               Provo, UT}
\newcommand{\yaff}{\normalsize Center for Demographic Studies,
                               Duke University, Durham, NC}

\begin{document}

\title{Grade of Membership Analysis:\\
    One Possible Approach to Foundations\thanks{
        This research was supported by grants
        from National Institute of Aging.}}

\author{
Mikhail Kovtun~~~~~~~
Igor Akushevich~~~~~~~
Kenneth G. Manton\\ \yaff
\and H. Dennis Tolley\\ \xaff
}

\markboth{M. Kovtun, I. Akushevich, K. G. Manton, and H. D. Tolley}
         {Grade of Membership Analysis: Approach to Foundations}

\maketitle

%-----------------------------------------------------------------------

\begin{abstract}
Grade of membership (GoM) analysis was introduced in 1974
\cite{Woodbury:1974} as a means of analyzing multivariate
categorical data.
Since then, it has been successfully applied to many problems.
The primary goal of GoM analysis is to derive
properties of individuals based on results of multivariate
measurements; such properties are given in the form of
the expectations of a hidden random variable (state of an individual)
conditional on the result of observations.

In this article, we present a new perspective for the GoM model,
based on considering distribution laws of observed random
variables as realizations of another random variable.
It happens that some moments of this new random variable
are directly estimable from observations.
Our approach allows us to establish a number of important
relations between estimable moments and values of interest,
which, in turn, provides a basis for a new numerical procedure.

\smallskip
\noindent
{\em {\bf Keywords: }
    Grade of membership analysis,
    latent structure analysis,
    multivariate categorical data,
    linear regression,
    multidimensional distribution.}

\smallskip
\noindent
{\em {\bf AMS 2000 subject classifications:}
    Primary 62H12; secondary 62J99.}
\end{abstract}

%=======================================================================

\section{Introduction}

The grade of membership (GoM) analysis was initially introduced
in \cite{Woodbury:1974};
the term ``Grade of Membership'' is due to this article.

GoM considers $J$ discrete measurements on each individual,
represented by random variables $X_1, \dots, X_J$,
with the set of outcomes of $j^\text{th}$ measurement
being $\{1, \dots, L_j\}$.

The goal of GoM analysis is to derive some properties
of an individual based on results of measurements.
We refer to this (general and informal) specification
of goals as {\em the General GoM Problem} (GGP.)
As GGP is a general concept, there may be many different
but reasonable answers to the problem.
The present article proposes one possible approach to GGP,
which leads to notable theoretical results and allows construction
of a novel numerical procedure.

Having GGP as its primary goal,
GoM differs from many other statistical
methods, whose goal is to discover some properties
of a population.
For example, in estimation of voting results
the most interesting fact is how many people will vote
for, or against (a candidate or an issue),
and it does not matter how a particular
individual votes.
In contrast, in making a medical diagnosis the health of
a particular individual is of interest, and it does not
matter (for a particular diagnosis) how prevalent
a particular health state is in a population.

Mathematically, one possibility to express GGP is to assume
that there exists a hidden continuous random variable $G$
representing knowledge about individuals derivable from observations
(in the diagnostic example, it is
the health state of an individual.)
Now one is interested in what might be said about value of $G$
based on observed values of $X_1, \dots, X_J$.
More specifically, values of interest are expectations of $G$
conditional on values of random variables $X_1, \dots, X_J$,
$\Expec(G \mid X_1=x_1, \dots, X_J=x_J)$.

Considering a continuous hidden random variables resembles
latent structure analysis in general, and
latent trait analysis in particular
(see \cite{Bartholomew:1999, Clogg:1995, Heinen:1996}.)
The connection between GoM and latent structure analysis
was mentioned in the literature (\cite{Haberman:1995}; see more details
in \cite{Erosheva:2002}.)
We prefer to keep the name ``grade of membership analysis''
because: (a) its primary goal differs from that
of the latent structure analysis, and (b) GoM uses a proprietary
technique and is based on special facts that are not
used in latent structure analysis.
However, we believe that techniques developed in the present article
and results obtained here might benefit the development
of latent structure analysis.

The main result of the present article, contained in
section \ref{sec:MainSystem},
is that the values of interest (i.e. conditional expectations and
conditional variances) are solution of system (\ref{eq:MainEqSys}),
and that under modest conditions, {\em only} values
of interest are solutions of this system.
Furthermore, as corollary \ref{cor:Main3} shows,
the system (\ref{eq:MainEqSys}) can be solved by two-step
process, every step of which consists of solving problem
of linear algebra.

GoM analysis (as well as many flavors of latent structure analysis)
employs an assumption that the problem under consideration
has lower dimensionality than observed data.
Our theorem \ref{th:Main2} and its corollary gives a way
to estimate this dimensionality directly
(which usually presents a substantual problem is such kind of analysis.)

An additional advantage of our approach is that it
not only establishes a way to estimate values of interest,
but also provides a ground for evaluation of confidence intervals
(not addressed in the present article.)

The rest of the article is organized as follows.

In section \ref{sec:TheProblem} we mathematically formulate
the problem and define related notions.
The central idea here (which is crucial for further results)
is to consider individual distribution laws
as realizations of another random variable, $\beta$.
We show that initial data are sufficient to estimate a set
of mixed moments of this distribution up to order $J$
(the number of measurements.)

In section \ref{sec:LowDimensionalDistributions} we consider the GoM
problem as a problem of finding a low-di\-men\-si\-o\-nal distribution
and obtain basic corollaries of this hypothesis.

In section \ref{sec:Regression} we consider a hypothesis that there
exists a linear regression of observed random variables $X_j$
on hidden random variable $G$.
We show that this hypothesis is essentially equivalent to
the one considered in previous section.

In section \ref{sec:RelationsBetweenMoments} we establish
relations between distributions and moments of $\beta$ and $G$,
and find transformation laws for changing their basis.
The main result of this section is equation (\ref{eq:MainEq}).

In section \ref{sec:MainSystem} we consider a system of equations
(\ref{eq:MainEqSys}). We show that values of interest are
always solutions of this system, and we establish
sufficient conditions under which the system (\ref{eq:MainEqSys})
has {\em only} such solutions.

In section \ref{sec:NumericalProcedure} we outline a numerical procedure
for estimating values of interest and discuss its properties.

%=======================================================================

\section{Preliminaries}
\label{sec:Preliminaries}

%-----------------------------------------------------------------------

\subsection{Notation}
\label{subsec:Notations}

$\mbZ$ is the set of integers, and $\mbR$ is the set of reals.
$\mbZ^+$ and $\mbR^+$ are subsets of positive,
and $\mbZ^{+0}$ and $\mbR^{+0}$ are subsets of nonnegative,
integers and reals, respectively.

For $m,n \in \mbZ$, $[m..n]$ denotes the set of integers
between $m$ and $n$:
$[m..n] = \{ z \in \mbZ \mid m \le z \le n \}$.
If $m > n$, $[m..n] = \varnothing$.

$\mbR^n$ is $n$-dimensional linear space over reals,
and $\mbS^n$ is a $(n-1)$-dimensional unit simplex in $\mbR^n$,
$\mbS^n = \{ x \in \mbR^n \mid x_i \ge 0 \text{ and } \sum_i x_i=1 \}$.

For a linear subspace $Q \subseteq \mbR^n$, $\dim(Q)$ denotes
its dimension.

For $x^1,\dots,x^p \in \mbR^n$, $\Lin(x^1,\dots,x^p)$ denotes
a linear subspace of $\mbR^n$ spanned by $x^1,\dots,x^p$,
and $\rng(x^1,\dots,x^p)$ denotes a rank of system of vectors
$x^1,\dots,x^p$ (thus,
$\rng(x^1,\dots,x^p) = \dim(\Lin(x^1,\dots,x^p))$.)

For $\alpha \in \mbR$ (or $\alpha \in \mbZ$) and $i \in [1..n]$,
$\boldsymbol{\alpha}_i$ denotes a vector from $\mbR^n$
($\mbZ^n$, respectively) with $i^\text{th}$ component equal $\alpha$
and all other components equal $0$. Dimensionality of
$\boldsymbol{\alpha}_i$ will be clear from context.

%-----------------------------------------------------------------------

\subsection{Support of measures}
\label{subsec:Support}

We consider only probabilistic measures defined on $\sigma$-algebra
of Borel sets of $\mbR^n$; a measure $\mu$ is a probabilistic measure,
if $\mu(\mbR^n) = 1$.

A {\em support} of measure $\mu$ is a closed set $A \subseteq \mbR^n$
such that $\mu(A) = 1$. We do not require a minimality of a support:
if $A$ is a support of $\mu$ and $A \subseteq A'$, $A'$ is closed,
then $A'$ also is a support of $\mu$.

We use $\Supp(\mu)$ to denote the set of all supports of $\mu$.
Thus, $A \in \Supp(\mu)$ means ``$A$ is a support of $\mu$.''
Note that $A \in \Supp(\mu)$ implies that $A$ is closed.

%-----------------------------------------------------------------------

\subsection{Indexing contingency tables and related objects}
\label{subsec:Indexing}

We need a way for indexing cells in a contingency table
and for other objects having similar structure.

A contingency table for a set of $J$ discrete measurements,
with $L_j$ possible outcomes for measurement $j$,
is a $J$-dimensional table having $L_j+1$ cells in dimension $j$.
Index for $j^\text{th}$ dimension ranges from $0$ to $L_j$.

More formally,
let $\mcL_\infty = \{ (l_1, \dots, l_J) \mid l_j \in \mbZ^+ \}$
and $\mcL_\infty^0 = \{ (l_1, \dots, l_J) \mid l_j \in \mbZ^{+0} \}$,
i.e. sets of $J$-dimensional vectors with positive
and, respectively, nonnegative integer components.
There is a one-to-one correspondence between sets of $J$ discrete
measurements and vectors in $\mcL_\infty$:
a vector $L = (L_1, \dots, L_J)$ describes a set of $J$ measurements,
in which measurement $j$ has $L_j$ outcomes.

For every $L \in \mcL_\infty$, let
$\mcL_L = \{ \ell \in \mcL_\infty \mid \ell_j \le L_j \}$
and $\mcL^0_L = \{ \ell \in \mcL^0_\infty \mid \ell_j \le L_j \}$.
If $L$ defines a set of measu\-rements, $\mcL_L$ is a set
of all possible outcomes of these measurements,
and $\mcL^0_L$ is a complete set of indices for
the contingency table.
In addition, for every $\smcJ \subseteq [1..J]$, let
$\mcL_L^{[\mcJ]} =
\{ \ell \in \mcL^0_L \mid \ell_j=0 \eqval j \in \smcJ \}$.
The set $\smcJ$ indicates measurements that we exclude from
consideration, and vector $\ell \in \mcL_L^{[\mcJ]}$ contains
results of all measurements except those listed in $\smcJ$.
Note that $\mcL_L^{[\varnothing]} = \mcL_L$ and 
$\mcL^0_L = \cup_{\mcJ \subseteq [1..J]} \mcL^{[\mcJ]}_L$.

Vector $\ell' \in \mcL_L^{[\mcJ]}$
may be considered as describing a family of outcomes
$\{ \ell \in \mcL_L \mid \ell_j=\ell'_j \text{ for } j \notin \smcJ \}$.
Abusing notation, we will also use $\ell'$ to denote this family.
More generally, we write $\ell' \in \ell''$ for
$\ell' \in \mcL_L^{[\mcJ']}$ and $\ell'' \in \mcL_L^{[\mcJ'']}$
whenever $\ell'_j = \ell''_j$ for all $j \notin \smcJ''$
(note $\ell' \in \ell''$ is possible only when
$\smcJ' \subseteq \smcJ''$.)
For $\ell \in \mcL_L^{[\mcJ]}$, let $\ell^{[\mcJ']}_{\phantom{0}}$
be a vector from $\mcL_L^{[\mcJ \cup \mcJ']}$ such that
$\ell^{[\mcJ']}_j = \ell_j$ for all $j \notin \smcJ'$.
We always have $\ell \in \ell^{[\mcJ']}_{\phantom{0}}$.
We write $\ell^{[j]}_{\phantom{0}}$, $\ell^{[j_1,j_2]}_{\phantom{0}}$,
etc. instead of $\ell^{[\{j\}]}_{\phantom{0}}$,
$\ell^{[\{j_1,j_2\}]}_{\phantom{0}}$, etc., respectively.

Let also set $|L| = \sum_j L_j$ and $|L^*| = \prod_j L_j$.

We always assume that the set of our measurements is
described by a vector $L$.
We drop index $L$ in notations $\mcL_L$ and $\mcL_L^{[\mcJ]}$
if it does not create an ambiguity.

A contingency table may be constructed for any sample
by putting in the cell with index $\ell$ the number of
individuals who (a) have outcome $\ell_j$ for measurement $j$
if $\ell_j \neq 0$; and (b) have arbitrary outcomes for all other
measurements. Let $N_\ell$ be a value in $\ell^\text{th}$ cell
of contingency table.
The usual summation rule for contingency tables in our notation
is: for any $\smcJ' \subseteq \smcJ \subseteq [1..J]$,
$\ell \in \smcJ$,
$N_\ell = \sum_{\ell'\in\mcJ \,:\, \ell'\in\ell} N_{\ell'}$.
Note that $N = N_{(0,\dots,0)}$ is the sample size.

A frequency table is obtained from a contingency table by dividing
the value in each cell by $N$. We use $f_\ell$ to denote a
value of $\ell^\text{th}$ cell of frequency table.
The above summation rule is applicable to frequency tables as well.

%=======================================================================

\section{The Problem}
\label{sec:TheProblem}

We consider a population of a potentially infinite number of individuals,
every individual being subject to $J$ measurements with discrete
outcomes. Without loss of generality, we may assume that outcomes
of $j^\text{th}$ measurement are $\{1, \dots, L_j\}$.

The results of measurements on individual $i$ is a random vector
$X^i = (X^i_1, \dots, X^i_J)$, which takes values in $\mcL_L$.
Such a random vector is described by a $|L|$-dimensional vector
of probabilities $\beta^i = (\beta^i_{jl})_{jl}$ ($j \in [1..J]$, and
for every $j$, $l \in [1..L_j]$), where $\beta^i_{jl} = \Prob(X^i_j=l)$.

These vectors of probabilities $\beta^i$ may themselves be considered
as realizations of a random vector $\beta$, with a distribution
described by probabilistic measure $\mu_\beta$ on $\mbR^{|L|}$.

We start with elementary properties, which may be directly
derived from definitions.

As $\beta_{jl}$ are probabilities, they satisfy

\begin{equation}
\label{eq:BetaCond}
\text{(a)} \quad \beta_{jl} \ge 0 \qquad \qquad
\text{(b)} \quad \text{for all $j$:}
                            \quad \sum_{l=1}^{L_j} \beta_{jl} = 1
\end{equation}

Thus, a product of simplices
$\mbS^L = \prod_j \mbS^{L_j} \subseteq \mbR^{|L|}$
is a support of the measure $\mu_\beta$,
$\mbS^L \in \Supp(\mu_\beta)$.

Together with random vectors $X^i$,
we consider a ``composite'' random vector
$X = (X_1, \dots, X_J)$:
on the first step, one randomly selects a vector of probabilities
$\beta$ (in accordance with measure $\mu_\beta$),
and on the second step, one randomly selects outcomes in
accordance with (selected on the first step) probabilities $\beta$.

According to our definitions, the conditional probability
for $X_j$ is:

\begin{equation}
\Prob \left( X_j = l \mid \beta \right) = \beta_{jl}
\end{equation}

\noindent
from which one obtains by the law of total probability

\begin{equation}
\label{eq:ModCond1}
\Prob \left( X_j = l \right) =
\int \Prob \left( X_j = l \mid \beta \right) \,\mu_\beta(d\beta) =
\int \beta_{jl} \,\mu_\beta(d\beta)
\end{equation}

We need more assumptions about $\mu_\beta$ to derive
useful properties of the model.
One reasonable assumption is ``local independence'':

\begin{enumerate}
\item[(G1)] Conditional on value of $\beta$,
random variables $X_1, \dots, X_J$ are mutually independent,
i.e. for every $\ell \in \mcL^0_{\phantom{L}}$

\begin{equation}
\Prob \Bigg( \bigwedge_{j \,:\, \ell_j \neq 0} X_j = \ell_j
                                                \Bigm| \beta \Bigg) =
\prod_{j \,:\, \ell_j \neq 0} \Prob \left(X_j=\ell_j \mid \beta \right)
\end{equation}
\end{enumerate}

A motivation for such assumption is that all ``randomness''
in $X^i_1, \dots, X^i_J$ comes from errors in measurements,
and error in one measurement does not depend on error in another
one.
Further, ``conditional on value of parameters'' means that we
are considering a group of individuals having the same values $\beta$;
thus, every individual in a group has the same vector of probabilities
$\beta$, and restriction of our random vector $X$ to this group
has the vector of probabilities $\beta$ as well;
as we assumed that for every individual random variables describing him
are independent, this should be true for a group
of identical (with respect to our random variables) individuals.
It is also wise to mention that the local independence assumption
is used in almost all variations of latent structure analysis.

With the independence assumption (G1), (\ref{eq:ModCond1}) may
be strengthened to:

\begin{equation}
\label{eq:ModCond2}
\forall \ell \in \mcL^0~ : \quad
\Prob \Bigg( \bigwedge_{j \,:\, \ell_j \neq 0} X_j = \ell_j \Bigg)
=
\int
    \Bigg( \prod_{j \,:\, \ell_j \neq 0} \beta_{j \ell_j} \Bigg)
\,\mu_\beta(d\beta)
\end{equation}

For every $\ell \in \mcL^0$, let the $\ell$-moment
of distribution $\mu_\beta$ be

\begin{equation}
\label{eq:MellDef}
M_\ell(\mu_\beta) =
    \int \Bigg( \prod_{j \,:\, \ell_j \neq 0} \beta_{j \ell_j} \Bigg)
\,\mu_\beta(d\beta)
\end{equation}

\noindent
In particular, we have
$\mcL_{\phantom{L}}^{[1,\dots,J]} = \{ (0,\dots,0) \}$, and
$M_{(0,\dots,0)}(\mu_\beta) = \int \mu_\beta(d\beta) = 1$.

Comparing (\ref{eq:MellDef}) with (\ref{eq:ModCond2}),
we see that the $\ell$-moment of distribution $\mu_\beta$
is equal to the probability of set of outcomes $\ell$.

For $\ell \in \mcL^{[\mcJ]}$, $M_\ell(\mu_\beta)$ is $J-|\smcJ|$ order
mixed moment of $\mu_\beta$. The set of $\ell$-moments for all
$\ell \in \mcL^0$ does not exhaust, however, the set of all moments
of order up to $J$ (for example, a moment
$\int \beta_{11} \beta_{12} \, \mu_\beta(d\beta)$
is not an $\ell$-moment.)
At the end of the section \ref{sec:MainSystem}
we shall discuss in more detail
whether $\{M_\ell(\mu_\beta)\}_\ell$ can determine
all moments of order up to $J$.

Basic statistical fact is that frequencies $f_\ell$ are consistent
and efficient estimators for $M_\ell(\mu_\beta)$.

The following proposition and its corollary is an equivalent
of the summation rule for contingency and frequency tables.

\begin{proposition}
\label{prop:sumMell}
Let $\smcJ' \subseteq \smcJ'' \subseteq [1..J]$.
Then for every $\ell'' \in \mcL^{[\mcJ'']}_{\phantom{L}}$

\begin{equation*}
M_{\ell''}(\mu_\beta) =
\sum_{\ell'\in\mcL^{[\mcJ']}_{\phantom{L}}~:~\ell'\in\ell''}
    M_{\ell'}(\mu_\beta)
\end{equation*}
\end{proposition}

\begin{prooff}
For every $j_0 \in \smcJ'' \setminus \smcJ'$
and for every $\ell' \in \mcL^{[\mcJ']}_{\phantom{L}}$ we have:

\allowdisplaybreaks{
\begin{multline*}
\sum_{\ell \in \mcL^{[\mcJ']}_{\phantom{L}}~:~
    \ell \in \ell'^{[j_0]}_{\phantom{0}}}
    M_{\ell}(\mu_\beta)
=
\sum_{\ell \in \mcL^{[\mcJ']}_{\phantom{L}}~:~
    \ell \in \ell'^{[j_0]}_{\phantom{0}}}
    \int \prod_{j \notin \mcJ'} \beta_{j \ell_j}
    \,\mu_\beta(d\beta)
=\\
\sum_{\ell \in \mcL^{[\mcJ']}_{\phantom{L}}~:~
    \ell \in \ell'^{[j_0]}_{\phantom{0}}}
    \int \beta_{j_0 \ell_{j_0}} \cdot
    \prod_{j \notin \mcJ'\cup\{j_0\}} \beta_{j \ell_j}
    \,\mu_\beta(d\beta)
=\\
\sum_{l = 1}^{L_{j_0}}
    \int \beta_{j_0 l} \cdot
    \prod_{j \notin \mcJ'\cup\{j_0\}} \beta_{j \ell_j}
    \,\mu_\beta(d\beta)
=\\
\int \left( \sum_{l = 1}^{L_{j_0}} \beta_{j_0 l} \right) \cdot
    \prod_{j \notin \mcJ'\cup\{j_0\}} \beta_{j \ell_j}
    \,\mu_\beta(d\beta)
=\\
\int 1 \cdot
    \prod_{j \notin \mcJ'\cup\{j_0\}} \beta_{j \ell_j}
    \,\mu_\beta(d\beta)
=
M_{\ell'^{[j_0]}}(\mu_\beta)
\end{multline*}
}

The rest of the proof is induction over the size of
$\smcJ'' \setminus \smcJ'$.
\end{prooff}

\begin{corollary}
\label{cor:sumMell}
For every $\smcJ \subseteq [1..J]$,
$\sum_{\ell \in \mcL^{[\mcJ]}_{\phantom{L}}} M_\ell(\mu_\beta) = 1$.
In particular,
$\sum_{\ell \in \mcL} M_\ell(\mu_\beta) = 1$.
\end{corollary}

Below we consider another two (essentially equivalent)
assumptions. The first one is that a support of $\mu_\beta$
is restricted to $(K-1)$-di\-men\-si\-o\-nal affine plane in $\mbR^{|L|}$.
The second assumption is that there exists a random variable $G$
taking values in $\mbR^K$ such that there exist a linear
regression of random variables $X_1,\dots,X_J$ on $G$.

%=======================================================================

\section{Low-dimensional distributions}
\label{sec:LowDimensionalDistributions}

The second assumption that we consider is:

\begin{enumerate}
\item[(G2$'$)]
    The support of $\mu_\beta$ is a $K$-dimensional linear
    subspace $Q$ of $\mbR^{|L|}$,
    and any proper subspace of $Q$ does not support
    $\mu_\beta$.
\end{enumerate}

We include the second clause (no proper subspace of $Q$
supports $\mu_\beta$) to avoid degenerate cases.
Any degenerate case may be considered as nondegenerate case
for some $K' < K$.

As $\mbS^L \in \Supp(\mu_\beta)$, the intersection
$P_\beta = Q \cap \mbS^L$ is necessarily nonempty,
and this intersection supports $\mu_\beta$.
In general, $P_\beta$ is $(K-1)$-dimensional polyhedral body,
which has at least $K$ vertices.
Let $\bar{P}_\beta$ be the $(K-1)$-dimensional affine space spanned
by $P_\beta$.

Let $\Lambda = \{ \lambda^1, \dots, \lambda^K \}$
be a linear basis of $Q$.
We also consider $\Lambda$ as a $|L| \times K$ matrix,

\begin{equation}
\Lambda =
\begin{pmatrix}
\lambda^1_{1 1}     & \dots     & \lambda^K_{1 1}   \\
                    & \vdots    &                   \\
\lambda^1_{J L_J}   & \dots     & \lambda^K_{J L_J} \\
\end{pmatrix}
\end{equation}

There exists considerable freedom in choosing $\Lambda$.
We shall exploit it by imposing constraints on $\Lambda$.
The first one is:

\begin{enumerate}
\item[$(\Lambda_0)$] For every $k$, $\lambda^k \in \bar{P}_\beta$.
\end{enumerate}

Let $g = (g_1, \dots, g_K)$ be a vector of coordinates of a point
$\beta \in Q$ in basis $\Lambda$,
i.e. $\beta = \sum_{k=1}^K g_k \lambda^k$.
Then $(\Lambda_0)$ implies

\begin{equation}
\sum_{k=1}^K g_k \lambda^k \in \bar{P}_\beta
    \quad \eqval \quad
\sum_{k=1}^K g_k = 1
\end{equation}

If $\Lambda$ and $\Lambda'$ are two bases of $Q$,
there exists a nondegenerate $K \times K$ matrix $A = (a^{k'}_k)_{k'k}$
such that $\Lambda' = \Lambda A$.

Using the fact that $(1,\dots,1)$ is left eigenvector of matrix $A$
corresponding to eigenvalue $1$
if and only if every column of $A$ sums to $1$, $\sum_k a^{k'}_k = 1$,
one easily obtains the following two propositions:

\begin{proposition}
Let both $\Lambda$ and $\Lambda'$ satisfy $(\Lambda_0)$.
Then $(1,\dots,1)$ is left eigenvector of matrix $A$
with eigenvalue $1$.
\end{proposition}

\begin{proposition}
Let $\Lambda$ satisfy $(\Lambda_0)$ and let $A$ be a nonsingular matrix
with left eigenvector $(1,\dots,1)$ with eigenvalue $1$.
Then $\Lambda' = \Lambda A$ satisfies $(\Lambda_0)$.
\end{proposition}

If $g$ is a coordinate vector of $\beta \in Q$ in basis $\Lambda$,
$\beta = \Lambda g$,
then the coordinate vector of $\beta$ in basis $\Lambda' = \Lambda A$
is $g' = A^{-1}g$.

\begin{remark}
In matrix expressions (like $\beta = \Lambda g$ above,) we always
assume that all vectors are columns.
$\blacksquare$
\end{remark}
\medskip

Every choice of a basis $\Lambda$ induces a linear map:

\begin{equation}
\label{eq:mapH}
H_\Lambda : \mbR^K \rightarrow Q, \qquad
H_\Lambda(g) = \sum_k g_k \lambda^k
\end{equation}

\noindent
Note that $\Lambda$ is a matrix of linear map $H_\Lambda$ with
respect to basis $\Lambda$ in $Q$ and standard unit
basis in $\mbR^K$.

When the basis $\Lambda$ satisfies $(\Lambda_0)$,
$H_\Lambda^{-1}(\bar{P}_\beta)$
is a unit affine plane $\bar{P}_g$ in $\mbR^K$,
$\bar{P}_g = \{ g \in \mbR^K \mid \sum_k g_k = 1 \}$,
and $P^\Lambda_g = H_\Lambda^{-1}(P_\beta)$ is
a convex $(K-1)$-dimensional polyhedron in $\bar{P}_g$.

The map $H_\Lambda$ allows us to introduce a measure $\mu^\Lambda_g$ on
$\bar{P}_g$, defined as:

\begin{equation}
\label{eq:muG}
\mu^\Lambda_g(B) = \mu_\beta(H_\Lambda(B)) \qquad
\text{for every Borel set } B \subseteq \bar{P}_g
\end{equation}

As $P_\beta \in \Supp(\mu_\beta)$, we have
$P^\Lambda_g \in \Supp(\mu^\Lambda_g)$.

Thus, we can replace integration over $P_\beta$ by integration
over $P^\Lambda_g$:

\begin{equation}
\int_{P_\beta} \phi(\beta) \mu_\beta(d\beta)  =
    \int_{P^\Lambda_g} \phi(H_\Lambda(g)) \mu^\Lambda_g(dg)
\end{equation}

\noindent
for every measurable function $\phi$.

\begin{remark}
We are trying to reflect in our notation all substantial
dependencies between objects.
Measure $\mu_\beta$ and polyhedron $P_\beta$, of course,
do not depend on the choice of $\Lambda$;
thus, no index $\Lambda$ in notation $\mu_\beta$ and $P_\beta$.
On the contrary, map $H$ defined by (\ref{eq:mapH})
(and consequently polyhedron $P_g$ and
measure $\mu_g$ defined by (\ref{eq:muG}))
substantially depends on the choice of $\Lambda$ --- so
we use notation $H_\Lambda$, $P^\Lambda_g$, and $\mu^\Lambda_g$.
However, we shall drop the index $\Lambda$ in the above notation
if it is obvious from the context.
$\blacksquare$
\end{remark}

%=======================================================================

\section{Linear regression hypothesis}
\label{sec:Regression}

A random variable $X_j$ has a finite range $[1..L_j]$,
on which no arithmetic operations are defined.
This prevents us from considering expectation, variance, etc.
of $X_j$.
To cope with this problem, we associate with every $X_j$
a random vector $Y_j$ taking values in $\mbR^{L_j}$
and defined as: if $X_j=l$, then $Y_j = \bsone_l$,
(recall that $\bsone_l$ is a $L_j$-dimensional vector with
$l^\text{th}$ component equals $1$, and all other components
equal $0$.)

There is an important connection between distributions of $X_j$
and $Y_j$: if $(\beta_{jl})_l$ is a vector of probabilities
of $X_j$, $\beta_{jl} = \Prob(X_j = l)$, then
$\Expec(Y_j) = (\beta_{jl})_l$
(here and below $\Expec(\cdot)$ denotes expectation.)
In general, for every condition $C$ we have
$\Expec(Y_j \mid C) = (\Prob(X_j=l \mid C)_l$.

\begin{remark}
As $Y_j$ is an $L_j$-dimensional vector,
$\Expec(Y_j)$ is also an $L_j$-dimensional vector.
We use $\Expec_m(\cdot)$ to denote $m^\text{th}$
component of vector expectation.
$\blacksquare$
\end{remark}
\medskip

Thus, we have

\begin{proposition}
\label{prop:ExpecProperties}
For every $j$ and for every condition $C$,
{\normalfont (a)} $\Expec_l(Y_j \mid C) \ge 0$,
and {\normalfont (b)} $\sum_l \Expec_l(Y_j \mid C) = 1$.
\end{proposition}

Now we can formulate an alternative form of assumption (G2):

\begin{enumerate}
\item[(G2$''$)]
    There exists a random vector $G$, defined on individuals
    and taking values in $\mbR^K$,
    such that:
    \begin{enumerate}
    \item There exists a joint distribution of $G$ and $X$.
    \item Local independence assumption holds, i.e. random variables
        $(X_1 \mid g)$, \dots, $(X_J \mid g)$ are mutually independent.
    \item For every $j$, a regression of $Y_j$ on $G$ is linear.
    \item For any $K' < K$ there is no random vector $G'$
        satisfying (a)--(c).
    \end{enumerate}
\end{enumerate}

Again, clause (d) is intended to prevent degenerate cases.

Clause (c) means that for every $j$, there exist
vectors $(\lambda^1_{jl})_l, \dots, (\lambda^K_{jl})_l$
such that

\begin{equation}
\label{eq:RegrHyp}
\Expec(Y_j \mid G=g) = \Big( \sum_k g_k \lambda^k_{jl} \Big)_l
\end{equation}

\noindent
or, in matrix form,

\begin{equation}
\Expec(Y_j \mid G=g) =
\begin{pmatrix}
\lambda^1_{j1}      & \dots     & \lambda^K_{j1}    \\
                    & \vdots    &                   \\
\lambda^1_{j L_j}   & \dots     & \lambda^K_{j L_j}
\end{pmatrix}
\cdot
\begin{pmatrix}
g_1     \\
\vdots  \\
g_K
\end{pmatrix} =
\Lambda_j \cdot g
\end{equation}

Taking into account the relation between $\Expec(Y_j)$ and
probability distribution of $X_j$, one obtains:

\begin{theorem}
(G2$'$) holds if, and only if, (G2$''$) holds.
\end{theorem}

The random vector $G$, if it exists, is not defined uniquely:
for every nondegenerate $K \times K$ matrix $A$,
random vector $G' = A^{-1}G$ also satisfies (G2$''$), as:

\begin{multline}
\Expec(Y_j \mid G'=g') =
\Expec(Y_j \mid AG'=Ag') =
\Expec(Y_j \mid G = Ag') =
\\
\Lambda_j \cdot (A \cdot g') =
(\Lambda_j \cdot A) \cdot g' =
\Lambda'_j \cdot g'
\end{multline}

This nonuniqueness corresponds to the nonuniqueness of the basis
for $Q$ discussed in section \ref{sec:LowDimensionalDistributions}.
Again, one may choose $G$ in such a way that $(\Lambda_0)$
is satisfied.

\begin{corollary}
In presence of $(\Lambda_0)$, the possible values of $G$
satisfy $\sum_k g_k = 1$.
In other words, $G$ takes values in a unit affine plane
$\bar{P}_g \subseteq \mbR^K$.
\end{corollary}

\begin{corollary}
In presence of $(\Lambda_0)$,
a set of possible values of $G$ is a bounded polyhedron
$P_g \subseteq \bar{P}_g$.
\end{corollary}

We are primarily interested in what can be said about
value of $G$ given outcomes of $X_1, \dots, X_J$.
The most interesting values are $\Expec(G \mid X = \ell)$
and $\Disper(G \mid X = \ell)$
(were $\Disper(\cdot)$ denotes variance.)
We shall derive equations for these values in the next section.
%\ref{sec:RelationsBetweenMoments}.

%=======================================================================

\section{Relations between $\mu_\beta$ and $\mu_g$}
\label{sec:RelationsBetweenMoments}

As (G2$'$) and (G2$''$) are equivalent, we refer to (either of)
them as (G2).

Under condition (G2) we have two distributions, $\mu_\beta$ and
$\mu^\Lambda_g$, connected by (\ref{eq:mapH}) and (\ref{eq:muG}).
In this section we establish further relations between 
$\mu_\beta$ and $\mu^\Lambda_g$.

Throughout this section, we assume that some basis $\Lambda$
of $Q$ is fixed.
We drop index $\Lambda$ in all notation;
however, the reader has to keep in mind that distribution $\mu_g$,
as well as all its moments, depend on $\Lambda$.

%-----------------------------------------------------------------------

\subsection{Unconditional moments}

We can express $\ell$-moments of $\mu_\beta$ via
moments of $\mu_g$.
Let $\smcJ \subseteq [1..J]$ and
$\ell \in \mcL^{[\mcJ]}_{\phantom{L}}$. Then:

\begin{multline}
\label{eq:MellMw}
M_\ell(\mu_\beta) =
%\\
\int_{P_\beta}
    \Bigg( \prod_{j \notin \mcJ} \beta_{j \ell_j} \Bigg)
\,\mu_\beta(d\beta) =
%\\
\int_{P_g}
    \Bigg( \prod_{j \notin \mcJ} \sum_k g_k \lambda^k_{j \ell_j} \Bigg)
\,\mu_g(dg) =
\\
\int_{P_g} \Bigg( \sum_{w \in \mcW^{[\mcJ]}_{\phantom{L}}}
    \bigg( \prod_{j \notin \mcJ} g_{w_j} \cdot
           \prod_{j \notin \mcJ} \lambda^{w_j}_{j \ell_j} \bigg)
\Bigg) \,\mu_g(dg) =
\\
\sum_{w \in \mcW^{[\mcJ]}_{\phantom{L}}}
    \Bigg(  
        \bigg(
            \int_{P_g} \prod_{j \notin \mcJ} g_{w_j}
            \,\mu_g(dg)
        \bigg)
        \cdot
        \prod_{j \notin \mcJ} \lambda^{w_j}_{j \ell_j}
    \Bigg) =
\\
\sum_{w \in \mcW^{[\mcJ]}_{\phantom{L}}}
    \bigg(
        M_w(\mu_g) \cdot
        \prod_{j \notin \mcJ} \lambda^{w_j}_{j \ell_j}
    \bigg)
\end{multline}

\noindent
Here
$\mcW^{[\mcJ]}_{\phantom{L}} =
\{ (w_1,\dots,w_J) \mid
    w_j \in [1..K] \text{ if } j \notin \smcJ,~
    w_j = 0 \text{ if } j \in \smcJ
\}$,
and for $w \in \mcW^{[\mcJ]}_{\phantom{L}}$,

\begin{equation}
M_w(\mu_g) =
\int_{P_g}
    \bigg( \prod_{j \,:\, w_j \neq 0} g_{w_j} \bigg)
\,\mu_g(dg)
\end{equation}

\noindent
is a $(J-|\smcJ|)^\text{th}$ order mixed moment of measure
$\mu^\Lambda_g$.

Note that $\mcW^{[\mcJ]}_{\phantom{L}} = \mcL^{[\mcJ]}_{(K,\dots,K)}$.
Thus, we freely apply to $\mcW$ all notations and conventions
developed for $\mcL$ in section \ref{subsec:Indexing}.

The sets of indices $\mcW^{[\mcJ]}_{\phantom{L}}$ are redundant
in the sense that different elements of $\mcW^{[\mcJ]}_{\phantom{L}}$
correspond to the same moments.
However, $\mcW^{[\mcJ]}_{\phantom{L}}$ has the following nice property:

\begin{proposition}
\label{prop:sumMw}
Let $\smcJ' \subseteq \smcJ'' \subseteq [1..J]$.
Then for every $w'' \in \mcW^{[\mcJ'']}_{\phantom{L}}$

\begin{equation*}
M_{w''}(\mu_g) =
\sum_{w' \in \mcW^{[\mcJ']}_{\phantom{L}}~:~w' \in w''}
    M_{w'}(\mu_g)
\end{equation*}
\end{proposition}

\begin{prooff}
Similar to the proof of proposition \ref{prop:sumMell}.
\end{prooff}

\begin{corollary}
\label{cor:sumMw}
For every $\smcJ \subseteq [1..J]$,
$\sum_{w \in \mcW^{[\mcJ]}_{\phantom{L}}} M_w(\mu_g) = 1$.
In particular,
$\sum_{w \in \mcW} M_w(\mu_g) = 1$.
\end{corollary}

To handle redundancy of $\mcW$, we introduce a new set
of indices.

Let $\mcV[J',K'] = \{ (v_1, \dots, v_{K'}) \mid v_k \in [0..J']
                                \text{ and } \sum_k v_k = J' \}$.
We write $\mcV[J']$ instead of $\mcV[J',K]$ and $\mcV$ instead
of $\mcV[J,K]$.

For every $\smcJ \subseteq [1..J]$ and for every
$v \in \mcV[J-|\smcJ|]$, let
$\mcW^{[\mcJ]}_v = \{ w \in \mcW^{[\mcJ]}_{\phantom{L}} \mid
\text{for} \linebreak[0]~
\text{every} \linebreak[0]~
k, \linebreak[0]~
w \linebreak[0]~
\text{contains} \linebreak[0]~
\text{exactly} \linebreak[0]~
v_k \linebreak[0]~
\text{components} \linebreak[0]~
\text{equal} \linebreak[0]~
k \}$.
Let also $C^{[\mcJ]}_v = |\mcW^{[\mcJ]}_v|$.

\begin{proposition}
\label{prop:valVC}
\begin{equation*}
\textnormal{(a)}\quad
|\mcV[J',K']| = \frac{(J'+K'-1)!}{J'! (K'-1)!}
\qquad \qquad
\textnormal{(b)}\quad
C^{[\mcJ]}_v = \frac{(J-|\smcJ|)!}{v_1! \dots v_K!}
\end{equation*}
\end{proposition}

\begin{prooff}
(a)
By induction over $J'+K'$ from a recurrent equality
$|\mcV[J',K']| = |\mcV[J'-1,K']| + |\mcV[J',K'-1]|$.

(b)
Let $J' = J - |\smcJ|$.
By direct computation one obtains:

\begin{equation*}
C^{[\mcJ]}_v = \binom{J'}{v_1} \cdot \binom{J'-v_1}{v_2}
    \dots \binom{J'-\sum_{k=1}^{K-1}v_k}{v_K}
\end{equation*}

\noindent
from which the statement of the proposition is straightforward.
\end{prooff}

One corollary to proposition \ref{prop:valVC} is that
$C^{[\mcJ]}_v$ depends on $\smcJ$ only through $|\smcJ|$,
and this value is contained in index $v$;
thus, we can safely drop index $[\smcJ]$ and write
just $C_v$.

As for every $w,w' \in \mcW^{[\mcJ]}_v$ we have
$M_w(\mu_g)=M_{w'}(\mu_g)$
for every measure $\mu_g$, we can define
{\em $v$-moments of a measure $\mu_g$} as

\begin{equation}
\label{eq:defMv}
M_v(\mu_g) =
M_{w_0}(\mu_g) =
\int \prod_k g_k^{v_k} \mu_g(dg)
\end{equation}

\noindent
and {\em normalized $v$-moments of a measure $\mu_g$} as

\begin{multline}
\label{eq:defnMv}
\bar{M}_v(\mu_g) = \sum_{w \in \mcW^{[\mcJ]}_v} M_w(\mu_g) =
C_v M_{w_0}(\mu_g) =
\\
C_v \int \prod_k g_k^{v_k} \mu_g(dg) =
C_v M_v(\mu_g)
\end{multline}

\noindent
In both equations, $w_0$ is an arbitrary element of $\mcW_v$.
Note that both $M_v(\mu_g)$ and $\bar{M}_v(\mu_g)$
do not depend on $\smcJ$.

$\mcV[J']$ is the smallest possible set of indices for $J'$-order mixed
moments of $\mu_g$. Multiplier $C_v$ in (\ref{eq:defnMv}) allows
us to obtain

\begin{proposition}
\label{prop:sumMv}
For every $\mu_g$,

\begin{equation*}
\sum_{v \in \mcV[J']} \bar{M}_v(\mu_g) =
\sum_{v \in \mcV[J']} C_v M_v(\mu_g) = 1
\end{equation*}
\end{proposition}

\begin{prooff}
Follows from proposition \ref{prop:sumMw}.
\end{prooff}

Now we can continue (\ref{eq:MellMw}):

\begin{multline}
\label{eq:MellMv}
\sum_{w \in \mcW^{[\mcJ]}_{\phantom{L}}}
    \bigg(
        M_w(\mu_g) \cdot
        \prod_{j \notin \mcJ} \lambda^{w_j}_{j \ell_j}
    \bigg) =
%\\
\sum_{v \in \mcV[J']}~ \sum_{w \in \mcW^{[\mcJ]}_v}
    \bigg(
        M_w(\mu_g) \cdot
        \prod_{j \notin \mcJ} \lambda^{w_j}_{j \ell_j}
    \bigg) =
\\
\sum_{v \in \mcV[J']}
    \bigg(
        M_v(\mu_g) \cdot
        \sum_{w \in \mcW^{[\mcJ]}_v}
            \prod_{j \notin \mcJ} \lambda^{w_j}_{j \ell_j}
    \bigg) =
%\\
\sum_{v \in \mcV[J']}
    \bigg(
        M_v(\mu_g) \cdot
        \sum_{w \in \mcW^{[\mcJ]}_v}
        \prod_{j \notin \mcJ} \lambda^{w_j}_{j \ell_j}
    \bigg) =
\\
\sum_{v \in \mcV[J']}
    \bigg(
        M_v(\mu^\Lambda_g) \cdot
        \Lambda(\smcJ,v,\ell)
    \bigg)
\end{multline}

\noindent
where $J' = J - |\smcJ|$ and

\begin{equation}
\label{eq:LambdaFunDef}
\Lambda(\smcJ,v,\ell) =
\sum_{w \in \mcW^{[\mcJ]}_v}
    \prod_{j \notin \mcJ} \lambda^{w_j}_{j \ell_j}
\end{equation}

%-----------------------------------------------------------------------

\subsection{Conditional moments}

For the joint distribution of $X=(X_1,\dots,X_J)$ and $G$
we have, on the one hand,

\begin{multline}
\label{eq:dProbG1}
d\Prob(G=g \wedge X=\ell) =
\Prob(X=\ell \mid G=g) \cdot d\Prob(G=g) =
\\
\Big( \prod_j \sum_k g_k \lambda^k_{j \ell_j} \Big) \,\mu_g(dg)
\end{multline}

\noindent
and, on the other hand,

\begin{multline}
\label{eq:dProbG2}
d\Prob(G=g \wedge X=\ell) =
d\Prob(G=g \mid X=\ell) \cdot \Prob(X=\ell) =
\\
d\Prob(G=g \mid X=\ell) \cdot M_\ell(\mu_\beta)
\end{multline}

Combining (\ref{eq:dProbG1}) and (\ref{eq:dProbG2}), one
obtains

\begin{equation}
\label{eq:ProbG3}
d\Prob(G=g \mid X=\ell) =
    \frac{\prod_j \sum_k g_k \lambda^k_{j \ell_j}}{M_\ell(\mu_\beta)}
    \,\mu_g(dg)
\end{equation}

Similarly, for every $\smcJ \subseteq [1..J]$ and for every
$\ell \in \mcL^{[\mcJ]}_{\phantom{L}}$,

\begin{equation}
\label{eq:ProbG4}
d\Prob(G=g \mid X = \ell) =
    \frac{\prod_{j\notin\mcJ} \sum_k g_k \lambda^k_{j \ell_j}}
         {M_\ell(\mu_\beta)}
    \,\mu_g(dg)
\end{equation}

\noindent
where for $\ell \in \mcL^{[\mcJ]}_{\phantom{L}}$,
$X = \ell$ means $\bigwedge_{j \not\in \mcJ} X_j = \ell_j$.

This allows us to conclude that the conditional distribution of
$G \mid X = \ell$ is absolutely continuous
with respect to measure $\mu_g$, and

\begin{equation}
\label{eq:CondProbDens}
p_\ell(g) = 
    \frac{\prod_{j\notin\mcJ} \sum_k g_k \lambda^k_{j \ell_j}}
         {M_\ell(\mu_\beta)}
\end{equation}

\noindent
is its probability density function.

Having this, we may write
(for every $J'$, $v \in \mcV[J']$, $\smcJ \subseteq [1..J]$,
$\ell \in \mcL^{[\mcJ]}_{\phantom{L}}$)
a $v$-order mixed moment of $G$ conditional on
$X=\ell$

\begin{equation}
\label{eq:MixMomG1}
\Expec(G^v \mid X=\ell) =
\int g^v p_\ell(g) \,\mu_g(dg)
\end{equation}

\noindent
where $g^v$ denotes $\prod_k g_k^{v_k}$.

A special case of equation (\ref{eq:MixMomG1}) for $v = (0,\dots,0)$ is

\begin{equation}
\label{eq:MixMomG2}
\Expec(G^{(0,\dots,0)} \mid X^{[\mcJ]}_{\phantom{L}}=\ell) =
\int g^{(0,\dots,0)} p_\ell(g) \,\mu_g(dg) =
\int p_\ell(g) \,\mu_g(dg) = 1
\end{equation}

Using equation (\ref{eq:MixMomG1}), we may obtain for
every $j \in \smcJ$ and every $l \in [1..L_j]$:

\begin{multline}
\label{eq:MainEq0}
\sum_k \lambda^k_{jl}
       \Expec(G^{v+\bsone_k} \mid X=\ell) =
\sum_k \lambda^k_{jl} \int g^{v+\bsone_k} p_\ell(g) \,\mu_g(dg) =
\\
\sum_k \int g^v (g_k \lambda^k_{jl}) p_\ell(g) \,\mu_g(dg) =
\int g^v \left( \sum_k g_k \lambda^k_{jl} \right)
    \frac{\prod_{j' \not\in \mcJ} \sum_k g_k \lambda^k_{j' \ell_{j'}}}
         {M_\ell(\mu_\beta)}
    \,\mu_g(dg) =
\\
\frac{M_{\ell+\bsl_j}(\mu_\beta)}{M_\ell(\mu_\beta)}
\int g^v \frac{\left( \sum_k g_k \lambda^k_{jl} \right)
           \prod_{j' \not\in \mcJ} \sum_k g_k \lambda^k_{j' \ell_{j'}}}
              {M_{\ell+\bsl_j}(\mu_\beta)}
    \,\mu_g(dg) =
\\
\frac{M_{\ell+\bsl_j}(\mu_\beta)}{M_\ell(\mu_\beta)}
\Expec(G^v \mid X=\ell+\bsl_j)
\end{multline}

By multiplying both sides of (\ref{eq:MainEq0}) by
$M_\ell(\mu_\beta)$ one obtains:

\begin{equation}
\label{eq:MainEq}
\sum_k \lambda^k_{jl} \cdot
\left(
    M_\ell(\mu_\beta) \cdot
    \Expec(G^{v+\bsone_k} \mid X=\ell)
\right) =
M_{\ell+\bsl_j}(\mu_\beta) \cdot \Expec(G^v \mid X=\ell+\bsl_j)
\end{equation}

Equation (\ref{eq:MainEq}) is the main fact that allows
us to establish a numerical procedure to estimate
conditional expectations.
This equation holds for every
$J' \ge 0$, $v \in \mcV[J']$, $\smcJ \subseteq [1..J]$,
$\ell \in \mcL^{[\mcJ]}_{\phantom{L}}$, $j \in \smcJ$,
and $l \in [1..L_j]$.

Although equation (\ref{eq:MainEq}) holds for every $J'$ and $\smcJ$,
the most important case is $J' + |\smcJ| < J$:
as we shall see in section \ref{sec:MainSystem},
only conditional moments $\Expec(G^v \mid X=\ell)$,
$v \in \mcV[J']$, $\ell \in \mcL^{[\mcJ]}_{\phantom{L}}$,
with $J' + |\smcJ| \le J$ may be identified from data.

%-----------------------------------------------------------------------

\subsection{Conditional variance}

To make use of conditional expectations, one would like
to know variance of $G$ conditional on outcomes of measurements.
It is not hard to express variance via conditional moments:
%(considered in previous section):

\begin{multline}
\label{eq:Disper}
\Disper_k(G \mid X=\ell) =
\int \left( g_k - \Expec_k(G \mid X=\ell) \right)^2 p_\ell(g)
\,\mu_g(dg) =
\\
\Expec(G^{\bstwo_k} \mid X=\ell) -
              \Expec^2(G^{\bsone_k} \mid X=\ell)
\end{multline}

As we shall show below, $\Expec(G^{\bstwo_k} \mid X=\ell)$
can be identified only for $\ell$ having at least two components
equal $0$;
thus, the same condition applies to identifiability of
$\Disper(G \mid X=\ell)$.

%-----------------------------------------------------------------------

\subsection{Change of basis}

Let $\Lambda' = \Lambda A$ be another basis of $Q$.
Here $A$ is nonsingular $K \times K$ matrix,

\begin{equation}
A = \begin{pmatrix}
        a^1_1 & \dots  & a^K_1  \\
              & \vdots &        \\
        a^1_K & \dots  & a^K_K
    \end{pmatrix},
\quad \text{and} \quad
A^{-1} = \begin{pmatrix}
        \bar{a}^1_1 & \dots  & \bar{a}^K_1  \\
                    & \vdots &              \\
        \bar{a}^1_K & \dots  & \bar{a}^K_K
    \end{pmatrix}
\end{equation}

As it was mentioned above, if a vector $\beta \in Q$ has
coordinates $g$ in basis $\Lambda$, then it has coordinates
$g' = A^{-1} g$ in basis $\Lambda'$.
Thus, $A^{-1}$ is a matrix of transition from coordinates $g$
to coordinates $g'$.
A question of interest is how the moments of $G$ are changed
under this transition.

We start with moments $M_w$ for $w \in \mcW$. Let $M'_w$ be
a moment calculated in coordinates $g'$. Then:

\allowdisplaybreaks{
\begin{multline}
\label{eq:Tensor1}
M'_w = M'_{(k_1,\dots,k_J)} =
\int g'_{k_1} \dots g'_{k_J} \,\mu'_g(dg') =
\\
\int (\bar{a}^1_{k_1} g_1 + \dots + \bar{a}^K_{k_1} g_K) \dots
     (\bar{a}^1_{k_J} g_1 + \dots + \bar{a}^K_{k_J} g_K) \,\mu_g(dg) =
\\
\sum_{m_1 = 1}^K \dots \sum_{m_J = 1}^K
    \bar{a}^{m_1}_{k_1} \dots \bar{a}^{m_J}_{k_J}
    \int g_{m_1} \dots g_{m_J} \,\mu_g(dg) =
\\
\sum_{m_1 = 1}^K \dots \sum_{m_J = 1}^K
    \bar{a}^{m_1}_{k_1} \dots \bar{a}^{m_J}_{k_J}
    M_{(m_1,\dots,m_J)}
\end{multline}
}

\noindent
which suggests that $\{M_w\}_{w \in \mcW}$ is a covariant tensor
of rank $J$. Employing Einstein's convention for summation,
(\ref{eq:Tensor1}) may be rewritten,

\begin{equation}
M'_{(k_1,\dots,k_J)} =
\bar{a}^{m_1}_{k_1} \dots \bar{a}^{m_J}_{k_J} M_{(m_1,\dots,m_J)}
\end{equation}

Tensor $\{M_w\}_{w \in \mcW}$ is symmetric, and
$\{M_v\}_{v \in \mcV}$ is a set of its essential components
(as for any $w \in v$, $M_w = M_v$.)
Transformation rules for $M_v$ have form:

\begin{equation}
\label{eq:Tensor2}
M'_v = \frac{1}{C_v} \sum_{v'}
    \left( \sum_{w \in v} \sum_{w' \in v'}
           \bar{a}^{w'_1}_{w_1} \dots \bar{a}^{w'_J}_{w_J}
    \right) M_{v'}
\end{equation}

For the general case of conditional moments of arbitrary order,
one obtains

\begin{equation}
\label{eq:Tensor3}
\Expec(G'^v \mid X=\ell) =
\frac{1}{C_v}
\sum_{v'}
    \left( \sum_{w \in v} \sum_{w' \in v'}
           \bar{a}^{w'_1}_{w_1} \dots \bar{a}^{w'_{J'}}_{w_{J'}}
    \right)
    \Expec(G^{v'} \mid X=\ell)
\end{equation}

\noindent
Here $v \in \mcV[J']$ for some $J'$, $v'$ ranges over $\mcV[J']$,
$w$ and $w'$ are restricted to the set
$\mcW[J'] = \{ (w_1,\dots,w_{J'}) \mid w_j \in [1..K] \}$,
and $w \in v$ means ``for every $k$, $w$ contains exactly $v_k$
components equal to $k$.''
%Of course, of the most interest are those cases of (\ref{eq:Tensor3}),
%which contain expectations that occur in equations (\ref{eq:MainEq}).

%=======================================================================

\section{Main system of equations}
\label{sec:MainSystem}

Consider a system of equations,

\begin{equation}
\label{eq:MainEqSys}
\begin{cases}
\sum_k \alpha^k_{jl} h^{v+\bsone_k}_\ell = h^v_{\ell+\bsl_j}, \quad
    & J' \in [0..J-1],~~
      v \in \mcV[J'],\\
    & \smcJ \subseteq [1..J] \,:\, |\smcJ|>J',~~
      \ell \in \mcL^{[\mcJ]}_{\phantom{L}},\\
    & j \in \smcJ,~~
      l \in [1..L_j]
\\[3pt]
h^{(0,\dots,0)}_\ell = M_\ell,
    & \ell \in \mcL^0
\\[3pt]
\sum_{v \in \mcV[J']} C_v h^v_{(0,\dots,0)} = 1,
    & J' \in [0..J]
\end{cases}
\end{equation}

\noindent
with respect to unknowns $\alpha^k_{jl}$
and $h^v_\ell$.

Equations (\ref{eq:MainEq}) and (\ref{eq:MixMomG2})
together with proposition \ref{prop:sumMv}
give us

\begin{theorem}
\label{th:Main1}
Let $\{M_\ell\}_{\ell \in \mcL^0}$ be a set of $\ell$-moments
of distribution $\mu_\beta$, which satisfies (G1) and (G2).
Let also $\{\lambda^k\}_k$, $\lambda^k = (\lambda^k_{jl})_{jl}$
be some basis of the support of $\mu_\beta$,
and $\Expec(G^v \mid X=\ell)$ be
conditional moments calculated with respect to this basis.

Then $\alpha^k_{jl} = \lambda^k_{jl}$ and 
$h^v_\ell = M_\ell \cdot \Expec(G^v \mid X=\ell)$
give a solution of system (\ref{eq:MainEqSys}).
\end{theorem}

In other words, all values we are interested in are solutions
of system (\ref{eq:MainEqSys}).
Below we establish sufficient conditions for the case
when (\ref{eq:MainEqSys}) has {\em only} such solutions.

For the sake of convenience, we (abusing language) shall speak
about ``solution $\alpha^1, \dots, \alpha^K$'', having in mind
``there exist $h^v_\ell$ such that $\alpha^1, \dots, \alpha^K$
together with $h^v_\ell$ compose a solution.''

Let $\alpha^1, \dots, \alpha^K$, $\alpha^k = (\alpha^k_{jl})_{jl}$,
and $h^v_\ell$ be a solution of (\ref{eq:MainEqSys}).
Let $\alpha'^1, \dots, \alpha'^K$ be any set of vectors such that
$\Lin(\alpha'^1, \dots, \alpha'^K) = \Lin(\alpha^1, \dots, \alpha^K)$.
In this case there exist a nonsingular $K \times K$ matrix
$A = (a^{k'}_k)_{k'k}$ such that
$(\alpha'^1, \dots, \alpha'^K) = (\alpha^1, \dots, \alpha^K) A$.
Let $A^{-1} = (\bar{a}^{k'}_k)_{k'k}$.
By straightforward computation one can show that
$\alpha'^1, \dots, \alpha'^K$ together with

\begin{equation*}
h'^v_\ell = 
\frac{1}{C_v}
\sum_{v'}
    \left( \sum_{w \in v} \sum_{w' \in v'}
           \bar{a}^{w'_1}_{w_1} \dots \bar{a}^{w'_{J'}}_{w_{J'}}
    \right)
    h^v_\ell
\end{equation*}

\noindent
also is a solution of (\ref{eq:MainEqSys}).

Thus, we can speak about {\em space of solutions}
$\Lin(\alpha^1, \dots, \alpha^K)$.
Note that at this point we have no arguments for
uniqueness of the space of solutions; moreover, we cannot even
claim that every space of solutions have the same dimension $K$.
In fact, in general case space of solutions is not unique.
However, in presence of sufficient conditions that we establish
below, the space of solutions is unique.

Consider equations from the first group of (\ref{eq:MainEqSys})
for $v = (0,\dots,0)$ and $\ell = (0,\dots,0)$:

\begin{equation}
\begin{cases}
\sum_k \alpha^k_{jl} h^{\bsone_k}_{(0,\dots,0)} =
    h^{(0,\dots,0)}_{\bsl_j}, \quad
& j \in [1..J],~~ l \in [1..L_j]
\end{cases}
\end{equation}

\noindent
and substitute values for $h^{(0,\dots,0)}_{\bsl_j}$
from the second group of (\ref{eq:MainEqSys}):

\begin{equation}
\label{eq:Basis1}
\begin{cases}
\sum_k \alpha^k_{jl} h^{\bsone_k}_{(0,\dots,0)} = M_{\bsl_j}, \quad
& j \in [1..J],~~ l \in [1..L_j]
\end{cases}
\end{equation}

As $h^{\bsone_k}_{(0,\dots,0)}$ do not depend on $j$ and $l$,
we obtain

\begin{proposition}
\label{prop:Basis1}
$(M_{\bsl_j})_{jl} \in \Lin(\alpha^1,\dots,\alpha^K)$
for every solution $\alpha^1, \dots, \alpha^K$ of (\ref{eq:MainEqSys}).
\end{proposition}

In other words, vector $(M_{\bsl_j})_{jl}$ belongs to every
space of solutions.

Applying similar considerations to the case $\ell = \bsl'_{j'}$
for some $j' \in [1..J]$, $l' \in [1..L_{j'}]$, we obtain:

\begin{equation}
\label{eq:Basis2}
\begin{cases}
\sum_k \alpha^k_{jl} h^{\bsone_k}_{\bsl'_{j'}} = M_{\bsl'_{j'}+\bsl_j},
\quad
& j \neq j',~~ l \in [1..L_j]
\end{cases}
\end{equation}

In system (\ref{eq:Basis2}) we have equations not for all $j,l$
but only for those in which $j \neq j'$.
Thus, (\ref{eq:Basis2}) does not give us a vector from a solution
space.
However, it allows us to claim that for every $j'$, $l'$,
a vector $(M_{\bsl'_{j'}+\bsl_j})_{jl \,:\, j \neq j'}$
(having $\sum_{j \neq j'} L_j$ components) may be extended
(by adding $L_{j'}$ components) to a $|L|$-dimensional vector
that belongs to $\Lin(\alpha^1,\dots,\alpha^K)$.

In general, for every $\ell \in \mcL^0 \setminus \mcL$ we have:

\begin{equation}
\label{eq:Basis3}
\begin{cases}
\sum_k \alpha^k_{jl} h^{\bsone_k}_{\ell} = M_{\ell+\bsl_j},
\quad
& \ell_j=0,~~ l \in [1..L_j]
\end{cases}
\end{equation}

\noindent
and thus we obtain further incomplete vectors
that may be completed to vectors belonging to
$\Lin(\alpha^1,\dots,\alpha^K)$.

Let us write vector $(M_{\bsl_j})_{jl}$ together with incomplete
vectors $(M_{\bsl'_{j'}+\bsl_j})_{jl \,:\, j \neq j'}$, etc.,
as columns of
a matrix, with places for which we do not have moments filled
by question marks. We refer to this incomplete matrix as to
{\em moment matrix}.
The moment matrix contains a column for every
$\ell \in \mcL^0 \setminus \mcL$.
Figure \ref{fig:1} gives an example of (part of) a moment matrix for
the case $J=3$, $L_1=L_2=L_3=2$.
Columns in this matrix correspond to
$\ell=(000)$, $(100)$, $(200)$, $(010)$, $(020)$, $(001)$, $(002)$,
$(110)$;
other columns are not shown.

\begin{figure}[ht]
\begin{equation*}
\begin{pmatrix}
M_{(100)}   &
    ?           &  ?            &
    M_{(110)}   &  M_{(120)}    &
    M_{(101)}   &  M_{(102)}    &
    ?           &  \cdots \phantom{\vdots}       \\
M_{(200)}   &
    ?           &  ?            &
    M_{(210)}   &  M_{(220)}    &
    M_{(201)}   &  M_{(202)}    &
    ?           &  \cdots \phantom{\vdots}       \\
M_{(010)}   &
    M_{(110)}   &  M_{(210)}    &
    ?           &  ?            &
    M_{(011)}   &  M_{(012)}    &
    ?           &  \cdots \phantom{\vdots}       \\
M_{(020)}   &
    M_{(120)}   &  M_{(220)}    &
    ?           &  ?            &
    M_{(021)}   &  M_{(022)}    &
    ?           &  \cdots \phantom{\vdots}       \\
M_{(001)}   &
    M_{(101)}   &  M_{(201)}    &
    M_{(011)}   &  M_{(021)}    &
    ?           &  ?            &
    M_{(111)}   &  \cdots \phantom{\vdots}       \\
M_{(002)}   &
    M_{(102)}   &  M_{(202)}    &
    M_{(012)}   &  M_{(022)}    &
    ?           &  ?            &
    M_{(112)}   &  \cdots \phantom{\vdots}      \\[-6pt]
&
\end{pmatrix}
\end{equation*}
\caption{\label{fig:1} Example of moment matrix}
\end{figure}

For a moment matrix $M$ let its completion $\bar{M}$ be a matrix
obtained from $M$ by replacing question marks by arbitrary numbers.
The above considerations give us

\begin{theorem}
\label{th:Main2}
Let distribution $\mu_\beta$ satisfy (G1) and (G2).
Then its moment matrix has a completion $\bar{M}$ such that
$\rng(\bar{M}) \le K$.
\end{theorem}

One may extend definition of rank to incomplete matrices
by setting it equal to the maximal size of nonzero minor,
which contains only known moments (i.e. does not contain question
marks.) It is easy to see that for every completion $\bar{M}$
of $M$, inequality $\rng(M) \le \rng(\bar{M})$ holds.
Thus,

\begin{corollary}
Let distribution $\mu_\beta$ satisfy (G1) and (G2).
Then $\rng(M) \le K$.
\end{corollary}

For $\mcK \subseteq \mcL^0 \setminus \mcL$,
let $M[\mcK]$ denote a matrix consisting of those columns
of moment matrix $M$ that correspond to elements of $\mcK$.

Now we are ready to formulate the third assumption regarding
distribution $\mu_\beta$:

\begin{enumerate}
\item[(G3)]
    There exist a subset of column indices
    $\mcK\subseteq \mcL^0 \setminus \mcL$ such that:
    \begin{enumerate}
    \item[(a)]
        For every two completions of moment matrix
        $\bar{M}'$ and $\bar{M}''$ satisfying
        $\rng(\bar{M}') \le K$ and $\rng(\bar{M}'') \le K$,
        the equality $\bar{M}'[\mcK] = \bar{M}''[\mcK]$ holds.
    \item[(b)]
        Let $\bar{M}$ be any completion of moment matrix 
        satisfying $\rng(\bar{M}) \le K$.
        Then $\rng(\bar{M}[\mcK]) = K$.
    \end{enumerate}
\end{enumerate}

Note that when (G3) holds, $\bar{M}[\mcK]$ is uniquely defined.

\begin{theorem}
\label{th:Main3}
Let distribution $\mu_\beta$ satisfy (G1), (G2), and (G3).
Then for every solution of system (\ref{eq:MainEqSys})
$\Lin(\alpha^1,\dots,\alpha^K) = \Lin(\bar{M}[\mcK])$
(where $\Lin(\bar{M}[\mcK])$ is a linear subspace of $\mbR^{|L|}$
spanned by columns of $\bar{M}[\mcK]$.)
\end{theorem}

\begin{prooff}
By theorem \ref{th:Main2}, for every solution of (\ref{eq:MainEqSys})
there exists a completion $\bar{M}'$ of $M$
such that $\Lin(\bar{M}') \subseteq \Lin(\alpha^1,\dots,\alpha^K)$.
Then $\rng(\bar{M}') = \dim(\Lin(\bar{M}')) \le
\dim(\Lin(\alpha^1,\dots,\alpha^K)) \le K$.
Thus, by (G3), $\bar{M}'[\mcK] = \bar{M}[\mcK]$, and consequently
$\Lin(\bar{M}[\mcK]) \subseteq \Lin(\bar{M}')$.
As $\dim(\Lin(\bar{M}[\mcK])) = \rng(\bar{M}[\mcK]) = K$, we obtain
$\Lin(\alpha^1,\dots,\alpha^K) = \Lin(\bar{M})$.
\end{prooff}

\begin{corollary}
\label{cor:Main3}
Let distribution $\mu_\beta$ satisfy (G1), (G2), and (G3).
Then:
\begin{enumerate}
\item[(a)]
    To obtain a solution of (\ref{eq:MainEqSys}), it is enough
    to take $\alpha^1, \dots, \alpha^K$ equal to any basis
    of $\bar{M}[\mcK]$ (e.g., equal to any $K$ linearly independent
    columns of $\bar{M}[\mcK]$.)
\item[(b)]
    Any other solution $\alpha'^1, \dots, \alpha'^K$ is obtained
    from the above one by multiplying it by nonsingular $K \times K$
    matrix.
\item[(c)]
    Every solution $\alpha^1, \dots, \alpha^K$ is a basis of
    $Q$, a support of $\mu_\beta$.
\end{enumerate}
\end{corollary}

\begin{prooff}
(a) and (b) are obvious.

To prove (c), consider that by theorem \ref{th:Main1}, every basis
of $Q$ is a solution of (\ref{eq:MainEqSys}). By (b), all solutions
are bases of the same linear subspace of $\mbR^{|L|}$.
Thus, every solution is basis of $Q$.
\end{prooff}

By theorem \ref{th:Main3} and its corollary, assumption (G3)
is sufficient to identify a support of $\mu_\beta$.
It looks like it is close to a necessary condition, as
in many cases where (G3) is violated, we were able to
construct a different distribution $\mu'_\beta$,
which has the same $\ell$-moments as $\mu_\beta$
(and therefore $\mu'_\beta$ is indistinguishable from $\mu_\beta$
based on available observations.)
However, the exact formulation of necessary conditions
for identifiability of support of $\mu_\beta$ is an open question.

To verify whether condition (G3) holds, it is enough to analyze
the moment matrix. Numerous practical methods might be suggested
to do such verification. Without going into details, we demonstrate
by example one possibility.

\begin{example}
\label{ex:1}
Consider a case $J=3$, $L_1=L_2=L_3=2$; thus, $\mbR^{|L|}=\mbR^6$.
Consider a distribution $\mu_\beta$ concentrated in three points,
$\beta^{(1)}$, $\beta^{(2)}$, and $\beta^{(3)}$, with every point having
probability $\frac{1}{3}$ (see figure \ref{fig:2}).
As $\beta^{(3)} = \frac{1}{2} \beta^{(1)} + \frac{1}{2} \beta^{(2)}$
and $\{ \beta^{(1)},\beta^{(2)}, \beta^{(3)} \} \in \Supp(\mu_\beta)$,
(G2) is satisfied for $K=2$.

\begin{figure}[ht!]
\begin{equation*}
\beta^{(1)} =
\begin{pmatrix}
    \beta^{(1)}_{11} \vphantom{\vdots} \\
    \beta^{(1)}_{12} \vphantom{\vdots} \\
    \beta^{(1)}_{21} \vphantom{\vdots} \\
    \beta^{(1)}_{22} \vphantom{\vdots} \\
    \beta^{(1)}_{31} \vphantom{\vdots} \\
    \beta^{(1)}_{32} \vphantom{\vdots} \\[-6pt]
    \phantom{.}
\end{pmatrix} =
\begin{pmatrix}
    1           \vphantom{\vdots} \\
    0           \vphantom{\vdots} \\
    \frac{1}{3} \vphantom{\vdots} \\
    \frac{2}{3} \vphantom{\vdots} \\
    \frac{4}{5} \vphantom{\vdots} \\
    \frac{1}{5} \vphantom{\vdots} \\[-6pt]
    \phantom{.}
\end{pmatrix},
~~
\beta^{(2)} =
\begin{pmatrix}
    \beta^{(2)}_{11} \vphantom{\vdots} \\
    \beta^{(2)}_{12} \vphantom{\vdots} \\
    \beta^{(2)}_{21} \vphantom{\vdots} \\
    \beta^{(2)}_{22} \vphantom{\vdots} \\
    \beta^{(2)}_{31} \vphantom{\vdots} \\
    \beta^{(2)}_{32} \vphantom{\vdots} \\[-6pt]
    \phantom{.}
\end{pmatrix} =
\begin{pmatrix}
    \frac{1}{9} \vphantom{\vdots} \\
    \frac{8}{9} \vphantom{\vdots} \\
    \frac{3}{5} \vphantom{\vdots} \\
    \frac{2}{5} \vphantom{\vdots} \\
    \frac{1}{4} \vphantom{\vdots} \\
    \frac{3}{4} \vphantom{\vdots} \\[-6pt]
    \phantom{.}
\end{pmatrix},
~~
\beta^{(3)} =
\begin{pmatrix}
    \beta^{(3)}_{11} \vphantom{\vdots} \\
    \beta^{(3)}_{12} \vphantom{\vdots} \\
    \beta^{(3)}_{21} \vphantom{\vdots} \\
    \beta^{(3)}_{22} \vphantom{\vdots} \\
    \beta^{(3)}_{31} \vphantom{\vdots} \\
    \beta^{(3)}_{32} \vphantom{\vdots} \\[-6pt]
    \phantom{.}
\end{pmatrix} =
\begin{pmatrix}
    \frac {5} {9} \vphantom{\vdots} \\
    \frac {4} {9} \vphantom{\vdots} \\
    \frac {7}{15} \vphantom{\vdots} \\
    \frac {8}{15} \vphantom{\vdots} \\
    \frac{21}{40} \vphantom{\vdots} \\
    \frac{19}{40} \vphantom{\vdots} \\[-6pt]
    \phantom{.}
\end{pmatrix}
\end{equation*}

\begin{equation*}
M =
\begin{pmatrix}
\frac{5}{9} &
    ?                &  ?                 &
    \frac {89} {405} &  \frac{136} {405}  &
    \frac{403}{1080} &  \frac{197}{1080}  &
    ?                &  \cdots \phantom{\vdots}       \\
\frac{4}{9} &
    ?                &  ?                 &
    \frac {20}  {81} &  \frac {16}  {81}  &
    \frac {41} {270} &  \frac {79} {270}  &
    ?                &  \cdots \phantom{\vdots}       \\
\frac{7}{15} &
    \frac {89} {405} &  \frac {20}  {81}  &
    ?                &  ?                 &
    \frac{397}{1800} &  \frac{443}{1800}  &
    ?                &  \cdots \phantom{\vdots}       \\
\frac{8}{15} &
    \frac{136} {405} &  \frac {16}  {81}  &
    ?                &  ?                 &
    \frac{137} {450} &  \frac{103} {450}  &
    ?                &  \cdots \phantom{\vdots}       \\
\frac{21}{40} &
    \frac{403}{1080} &  \frac {41} {270}  &
    \frac{397}{1800} &  \frac{137} {450}  &
    ?                &  ? [\leftarrow \frac{191}{960}]  &
    \frac{151}{1080} &  \cdots \phantom{\vdots}       \\
\frac{19}{40} &
    \frac{197}{1080} &  \frac {79} {270}  &
    \frac{443}{1800} &  \frac{103} {450}  &
    ?                &  ? [\leftarrow \frac{53}{192}] &
    \frac{259}{3240} &  \cdots \phantom{\vdots}      \\[-6pt]
&
\end{pmatrix}
\end{equation*}

\begin{equation*}
\alpha^1 =
\begin{pmatrix}
    \alpha^1_{11} \vphantom{\vdots} \\
    \alpha^1_{12} \vphantom{\vdots} \\
    \alpha^1_{21} \vphantom{\vdots} \\
    \alpha^1_{22} \vphantom{\vdots} \\
    \alpha^1_{31} \vphantom{\vdots} \\
    \alpha^1_{32} \vphantom{\vdots} \\[-6pt]
    \phantom{.}
\end{pmatrix} =
\begin{pmatrix}
    \frac {5} {9} \vphantom{\vdots} \\
    \frac {4} {9} \vphantom{\vdots} \\
    \frac {7}{15} \vphantom{\vdots} \\
    \frac {8}{15} \vphantom{\vdots} \\
    \frac{21}{40} \vphantom{\vdots} \\
    \frac{19}{40} \vphantom{\vdots} \\[-6pt]
    \phantom{.}
\end{pmatrix},
\qquad
\alpha^2 =
\begin{pmatrix}
    \alpha^2_{11} \vphantom{\vdots} \\
    \alpha^2_{12} \vphantom{\vdots} \\
    \alpha^2_{21} \vphantom{\vdots} \\
    \alpha^2_{22} \vphantom{\vdots} \\
    \alpha^2_{31} \vphantom{\vdots} \\
    \alpha^2_{32} \vphantom{\vdots} \\[-6pt]
    \phantom{.}
\end{pmatrix} =
\begin{pmatrix}
    \frac{197}{513} \vphantom{\vdots} \\
    \frac{316}{513} \vphantom{\vdots} \\
    \frac{443}{855} \vphantom{\vdots} \\
    \frac{412}{855} \vphantom{\vdots} \\
    \frac{191}{456} \vphantom{\vdots} \\
    \frac{265}{456} \vphantom{\vdots} \\[-6pt]
    \phantom{.}
\end{pmatrix}
\end{equation*}
\caption{\label{fig:2} Illustration to example \ref{ex:1}}
\end{figure}

The moment matrix $M$ of this distribution 
(which corresponds to moment matrix on figure \ref{fig:1})
is shown on figure \ref{fig:2}.

A submatrix of $M$ consisting of rows 3 and 4 and columns 1 and 2 is
nonsingular, and therefore $x$ and $y$ such that

\begin{equation*}
\text{column}1 \cdot x ~+~ \text{column}2 \cdot y ~=~ \text{column}7
\end{equation*}

\noindent
are uniquely defined;
they are $x=\frac{131}{160}$ and $y=-\frac{99}{160}$.
This allows construction of the only possible completion of column 7,
which is shown on figure \ref{fig:2}.

Thus, column 1 and (completed) column 7 give a basis for a support
of $\mu_\beta$. It is easy to see that
$\Lin(\text{column}1,\text{column}7) =
 \Lin(\beta^{(1)},\beta^{(2)},\beta^{(3)})$,
as one would expect.

Vectors ``$\text{column}1$'' and ``$\text{completed column}7$''
do not satisfy condition $(\Lambda_0)$.
To obtain a basis satisfying $(\Lambda_0)$,
one can take $\alpha^1 = \text{column}1$ and
$\alpha^2 = \text{column}7 \cdot \frac{40}{19}$.
Vectors $\alpha^1$ and $\alpha^2$ are shown on figure \ref{fig:2}.

(Calculations for this and subsequent examples were
done with Waterloo Maple$^{\text{TM}}$ v.7.00.)
$\blacksquare$
\end{example}

The second question is whether $h^v_\ell$ may be uniquely determined
from (\ref{eq:MainEqSys}) given a solution
$\alpha^1, \dots, \alpha^K$.
In general, the answer is negative: not all $h^v_\ell$ may be
uniquely determined. However, a number of the most important values
always may be determined uniquely, as the following theorem shows.

\begin{theorem}
\label{th:Main4}
Let $\alpha^1, \dots, \alpha^K$ be a solution of (\ref{eq:MainEqSys}),
and let set of index pairs $j_1 l_1, \dots, j_K l_K$,
with $l_k \in [1..L_{j_k}]$, be chosen so that the matrix
$(\alpha^{k'}_{j_k l_k})_{k'k}$ is nonsingular
(this is always possible as $\rng(\alpha^1, \dots, \alpha^K) = K$.)
Let $\smcJ_0 = \{ j_1, \dots, j_K \}$ (note that $|\smcJ_0|$ may be
less than $K$.)
Then:

\begin{enumerate}
\item[(a)]
    For every $\smcJ$ such that $\smcJ_0 \subseteq \smcJ$,
    for every $\ell \in \mcL^{[\mcJ]}_{\phantom{L}}$,
    and for every $k \in [1..K]$,
    the conditional expectation $\Expec(G_k \mid X = \ell)$
    is uniquely defined.
\item[(b)]
    Let, in addition, there exist $j_0 \not\in \smcJ_0$ and
    $l_0 \in [1..L_{j_0}]$ such that every $K \times K$ submatrix
    of $(K+1) \times K$ matrix
    $(\alpha^{k'}_{j_k l_k})_{k' \in [1..K], k \in [0..K]}$
    is nonsingular.
    Then for every $\smcJ$ such that
    $\smcJ_0 \cup \{j_0\} \subseteq \smcJ$,
    for every $\ell \in \mcL^{[\mcJ]}_{\phantom{L}}$,
    and for every $k \in [1..K]$,
    the conditional variance $\Disper(G_k \mid X = \ell)$
    is uniquely defined.
\end{enumerate}
\end{theorem}

\begin{prooff}
(a)
Consider a subsystem of (\ref{eq:MainEqSys}):

\begin{equation*}
\begin{cases}
\sum_{k'} \alpha^{k'}_{j_k l_k} h^{\bsone_{k'}}_\ell =
    M_{\ell + (\bsl_k)_{(j_k)}}, \quad
& k = 1,\dots,K
\end{cases}
\end{equation*}

\noindent
By theorem \ref{th:Main1},
$h^{\bsone_{k'}}_\ell =M_\ell \cdot \Expec(G_k \mid X = \ell)$
is a solution of this system, and by assumption of the theorem,
there are no other solutions.

(b)
By part (a) of the theorem, for every $k_0 \in [1..K]$
and every $k \in [1..K]$, values
$h^{\bsone_{k_0}}_{\ell + (\bsl_k)_{(j_k)}}$ are uniquely determined
from (\ref{eq:MainEqSys}).
Now consider a subsystem of (\ref{eq:MainEqSys}):

\begin{equation*}
\begin{cases}
\sum_{k'} \alpha^{k'}_{j_k l_k} h^{\bsone_{k_0} + \bsone_{k'}}_\ell =
    h^{\bsone_{k_0}}_{\ell + (\bsl_k)_{(j_k)}}, \quad
& k = 1,\dots,K;~~ k_0 = 1,\dots,K
\end{cases}
\end{equation*}

\noindent
By theorem \ref{th:Main1},
$h^{\bsone_{k_0}+\bsone_{k'}}_\ell =
M_\ell \cdot \Expec(G^{\bsone_{k_0}+\bsone_{k'}} \mid X = \ell)$
is a solution of this system, and by assumption of the theorem,
there are no other solutions. This is enough to calculate
$\Disper(G_k \mid X = \ell)$ using formula (\ref{eq:Disper}).
\end{prooff}

\begin{example}
\label{ex:2}
We continue example \ref{ex:1}.
Consider a subsystem of (\ref{eq:MainEqSys}):

\begin{equation*}
\begin{cases}
\alpha^1_{21} h^{(1,0)}_{(1,0,0)} + \alpha^2_{21} h^{(0,1)}_{(1,0,0)}
= M_{(1,1,0)}
\\[6pt]
\alpha^1_{22} h^{(1,0)}_{(1,0,0)} + \alpha^2_{22} h^{(0,1)}_{(1,0,0)}
= M_{(1,2,0)}
\end{cases},
\quad \text{or} \quad
\begin{cases}
\frac{7}{15} h^{(1,0)}_{(1,0,0)} + \frac{443}{855} h^{(0,1)}_{(1,0,0)}
= \frac{89}{405}
\\[6pt]
\frac{8}{15} h^{(1,0)}_{(1,0,0)} + \frac{412}{855} h^{(0,1)}_{(1,0,0)}
= \frac{136}{405}
\end{cases}
\end{equation*}

Solving this system gives

\begin{equation*}
h^{(1,0)}_{(1,0,0)} = \frac{131}{99}, \qquad
h^{(0,1)}_{(1,0,0)} = - \frac{76}{99}
\end{equation*}

\noindent
and, as $h^v_\ell = M_\ell \cdot \Expec(G^v \mid X = \ell)$,

\begin{equation*}
\Expec(G^{(1,0)} \mid X = (1,0,0)) = \frac{131}{55}, \qquad
\Expec(G^{(0,1)} \mid X = (1,0,0)) = - \frac{76}{55}
\end{equation*}

Considering similar subsystems, one obtains, in particular,

\begin{align*}
h^{(1,0)}_{(1,0,1)} =   \frac{3089}{2970}, \quad
h^{(0,1)}_{(1,0,1)} = - \frac{2641}{3960}, \quad
h^{(1,0)}_{(1,0,2)} =   \frac {841}{2970}, \quad
h^{(0,1)}_{(1,0,2)} = - \frac {133}{3960}
\end{align*}

Substituting these values into subsystems,

\begin{equation*}
\begin{cases}
\alpha^1_{31} h^{(2,0)}_{(1,0,0)} +
\alpha^2_{31} h^{(1,1)}_{(1,0,0)} = h^{(1,0)}_{(1,0,1)}
\\[6pt]
\alpha^1_{32} h^{(2,0)}_{(1,0,0)} +
\alpha^2_{32} h^{(1,1)}_{(1,0,0)} = h^{(1,0)}_{(1,0,2)}
\end{cases},
\quad
\begin{cases}
\alpha^1_{31} h^{(1,1)}_{(1,0,0)} +
\alpha^2_{31} h^{(0,2)}_{(1,0,0)} = h^{(0,1)}_{(1,0,1)}
\\[6pt]
\alpha^1_{32} h^{(1,1)}_{(1,0,0)} +
\alpha^2_{32} h^{(0,2)}_{(1,0,0)} = h^{(0,1)}_{(1,0,2)}
\end{cases}
\end{equation*}

\noindent
one finds,

\begin{equation*}
h^{(2,0)}_{(1,0,0)} =   \frac{3323} {726}, \qquad
h^{(1,1)}_{(1,0,0)} = - \frac{7087}{2178}, \qquad
h^{(0,2)}_{(1,0,0)} =   \frac{1895} {726}
\end{equation*}

\noindent
and thus,

\begin{equation*}
\Expec( G^{(2,0)} \mid X = (1,0,0) ) = \frac{9969}{1210}, \qquad
\Expec( G^{(0,2)} \mid X = (1,0,0) ) = \frac{1083} {242}
\end{equation*}

This allows us calculate conditional variances (using formula
(\ref{eq:Disper})):

\begin{equation*}
\Disper( G_1 \mid X = (1,0,0) ) = \frac{15523}{6050}, \qquad
\Disper( G_2 \mid X = (1,0,0) ) = \frac{15523}{6050}
\end{equation*}

Table \ref{tab:1} summarize conditional expectations and conditional
variances that may be calculated in our example.
Although all values are {\em exact} rational numbers, we used
decimal notation to make comparison of values easier.
We also put standard deviations in the table instead of variances.

As we have mentioned, there are many choices for 
basis for the support of distribution $\mu_\beta$.
Another possibility is to take $\{ \beta^{(1)}, \beta^{(2)} \}$
as a basis. The result of calculations in this basis is given
in table \ref{tab:2}. One can see that, although numbers are different,
their relative position remains the same.
$\blacksquare$
\end{example}

\begin{table}[ht]
\caption{\label{tab:1}
    Conditional expectations and standard deviations
    calculated in basis $\{ \alpha^1, \alpha^2 \}$.
}
\begin{tabular}{|c||c|c||c|c|}
\hline
\vphantom{$\Big(\Big)$}
$\ell$ & $\Expec(G_1 \mid X=\ell)$ & $\sigma(G_1 \mid X=\ell)$
       & $\Expec(G_2 \mid X=\ell)$ & $\sigma(G_2 \mid X=\ell)$
\\
\hline\hline
\vphantom{$\Big(\Big)$}
$(1,0,0)$ & 2.3818 & 1.6018 & -1.3818 & 1.6018
\\
\hline
\vphantom{$\Big(\Big)$}
$(2,0,0)$ & -0.7273 & 1.2214 & 1.7273 & 1.2214
\\
\hline
\vphantom{$\Big(\Big)$}
$(0,1,0)$ & 0.5065 & 2.0571 & 0.4935 & 2.0571
\\
\hline
\vphantom{$\Big(\Big)$}
$(0,2,0)$ & 1.4318 & 2.0709 & -0.4318 & 2.0709
\\
\hline
\vphantom{$\Big(\Big)$}
$(0,0,1)$ & 1.9048 & 1.9122 & -0.9048 & 1.9122
\\
\hline
\vphantom{$\Big(\Big)$}
$(0,0,2)$ & 0.0000 & 1.8642 & 1.0000 & 1.8642
\\
\hline
\end{tabular}
\medskip
\end{table}

\begin{table}[ht]
\caption{\label{tab:2}
    Conditional expectations and standard deviations
    calculated in basis $\{ \beta^{(1)}, \beta^{(2)} \}$.
}
\begin{tabular}{|c||c|c||c|c|}
\hline
\vphantom{$\Big(\Big)$}
$\ell$ & $\Expec(G_1 \mid X=\ell)$ & $\sigma(G_1 \mid X=\ell)$
       & $\Expec(G_2 \mid X=\ell)$ & $\sigma(G_2 \mid X=\ell)$
\\
\hline\hline
\vphantom{$\Big(\Big)$}
$(1,0,0)$ & 0.7667 & 0.3091 & 0.2333 & 0.3091
\\
\hline
\vphantom{$\Big(\Big)$}
$(2,0,0)$ & 0.1667 & 0.2357 & 0.8333 & 0.2357
\\
\hline
\vphantom{$\Big(\Big)$}
$(0,1,0)$ & 0.4048 & 0.3970 & 0.5952 & 0.3970
\\
\hline
\vphantom{$\Big(\Big)$}
$(0,2,0)$ & 0.5833 & 0.3997 & 0.4167 & 0.3997
\\
\hline
\vphantom{$\Big(\Big)$}
$(0,0,1)$ & 0.6746 & 0.3690 & 0.3254 & 0.3690
\\
\hline
\vphantom{$\Big(\Big)$}
$(0,0,2)$ & 0.3070 & 0.3598 & 0.6930 & 0.3598
\\
\hline
\end{tabular}
\medskip
\end{table}

\begin{remark}
The standard deviations in the above example
are relatively large. This is direct consequence of the fact
that in this example we have too small number of measurements.
When number of measuremnents increases, the standard deviation
becomes smaller and smaller.
$\blacksquare$
\end{remark}
%\medskip

\begin{remark}
Theorem \ref{th:Main4} guarantees that it is always possible to
find $J-K$ measurements such that expectations of $G$ conditional
on outcomes of these measurements may be uniquely determined
from the system (\ref{eq:MainEqSys}).
The possibility of determining conditional variances is not
guaranteed by this theorem, however.
In many practical cases that we have investigated, conditions of the
part (b) of theorem \ref{th:Main4} are satisfied, and conditional
variances can be found (as in example \ref{ex:2}.)
The exact conditions for determinability
of conditional variances is an open question.
$\blacksquare$
\end{remark}
%\medskip

\begin{remark}
By computations similar to used in (\ref{eq:MellMw}) and
(\ref{eq:MellMv}), one obtains for every family of $J' \le J$
index pairs $j_1 l_1, \dots, j_{J'} l_{J'}$ with $l_p \in [1..L_{j_p}]$
($j_p$ is not necessarily different from $j_{p'}$ for $p \neq p'$)

\begin{equation*}
\int \beta_{j_1 l_1} \cdot \ldots \cdot \beta_{j_{J'} l_{J'}}
    \,\mu_\beta(d\beta) =
\sum_{v \in \mcV[J']} M_v(\mu_g) \cdot
    \tilde{\Lambda}(v,j_1,l_1,\dots,j_{J'},l_{J'})
\end{equation*}

\noindent
where $\tilde{\Lambda}(v,j_1,l_1,\dots,j_{J'},l_{J'})$ depends only
on $\lambda^k_{jl}$.
Thus, if the system (\ref{eq:MainEqSys}) allows unique determination of
all unknowns $h^v_\ell$, all moments of order up to $J$ of $\mu_\beta$
can be identified.
This is the case, for instance, in the example \ref{ex:2}.

We do not know now whether there exist some regular conditions
under which the system (\ref{eq:MainEqSys}) has a unique solution
(modulo change of basis.)
Examples that we have considered suggest that in a regular case
system (\ref{eq:MainEqSys}) {\em never} has a unique solution
whenever $K > L_j$ at least for one $j$.
(However, as theorem \ref{th:Main4} shows, many values of interest
always may be uniquely determined.)
The exact description of parameters that may be uniquely identified
based on system (\ref{eq:MainEqSys}), and to what degree the freedom
in choosing other parameters may be reduced,
is a subject for further investigation.
$\blacksquare$
\end{remark}

%=======================================================================

\section{Numerical procedure}
\label{sec:NumericalProcedure}

We have established a number of
{\em precise} relations between values of interest
(i.e. expectation and variance of hidden random vector $G$
conditional on outcomes of measurements)
and moments of (unknown) distribution $\mu_\beta$, which are directly
estimable from observations.
The most important of these relations are given by equations
(\ref{eq:MainEq}), and by system of equations (\ref{eq:MainEqSys}).
This relations suggest a numerical procedure for estimation of
values of interest.

As was mentioned above, sample frequences $f_\ell$ are
consistent estimators for moments $M_\ell(\mu_\beta)$.
Thus, applying the least squares method to the system

\begin{equation}
\label{eq:MainEqSysFreq}
\begin{cases}
\sum_k \alpha^k_{jl} h^{v+\bsone_k}_\ell = h^v_{\ell+\bsl_j}, \quad
    & J' \in [0..J-1],~~
      v \in \mcV[J'],\\
    & \smcJ \subseteq [1..J] \,:\, |\smcJ|>J',~~
      \ell \in \mcL^{[\mcJ]}_{\phantom{L}},\\
    & j \in \smcJ,~~
      l \in [1..L_j]
\\[3pt]
h^{(0,\dots,0)}_\ell = f_\ell,
    & \ell \in \mcL^0
\\[3pt]
\sum_{v \in \mcV[J']} C_v h^v_{(0,\dots,0)} = 1,
    & J' \in [0..J]
\end{cases}
\end{equation}

\noindent
one obtains consistent estimators for a basis $\{\lambda^k\}_k$
and conditonal expectations of $G$.

The consistency of estimators obtained from (\ref{eq:MainEqSysFreq})
is almost straightforward corollary to consistency of estimators
$f_\ell$. The rate of convergence is more delicate question
(as a rate of convergence of $f_\ell$ depends on $\ell$,)
and deserves separate investigation.

Theorem \ref{th:Main3} suggests another, two-step way for finding
solutions of (\ref{eq:MainEqSysFreq}).
On the first step, one finds a basis from frequency matrix
(i.e. moment matrix with frequences substituted for moments.)
After basis is obtained, (\ref{eq:MainEqSysFreq}) turns to be
a linear system with respect to $h^v_\ell$.
This way requires significantly less computations, but
its convergence properties have to be more carefully investigated.

One question regarding numerical procedure is the choice of
value of $K$ for which system (\ref{eq:MainEqSysFreq}) should
be solved. Theorem \ref{th:Main2} and its corollary suggest
that one has to take $K$ equal to the rank of the frequency
matrix (modulo possible deviations of frequencies from the
true moments.)

Another question is how a numerical algorithm has to deal with is
nonuniqueness of basis $\{\lambda^k\}_k$.
In general, there are $K^2$ degrees of freedom in choice of a basis.
Imposing condition $(\Lambda_0)$ reduces this number to $K(K-1)$.
One can consider additional restrictions on choice of basis:

\begin{enumerate}
\item[$(\Lambda_1)$]
    For every $k$, unconditional expectation $\Expec_k(G)$ equals
    $\frac{1}{K}$.
\item[$(\Lambda_2)$]
    The map $H_\Lambda$ is isometry of $\bar{P}_g$ and $\bar{P}_\beta$
    (with respect to euclidean distance.)
\end{enumerate}

The firts one corresponds to restricting transformations of $\bar{P}_g$,
described by matrix $A$ (introduced in section
\ref{sec:LowDimensionalDistributions}) to those having
the ``center'' point $(\frac{1}{K},\dots,\frac{1}{K})$ of $\bar{P}_g$
fixed.
The second restriction guarantees that variances do not depend
on the choice of basis,
and variances calculated in $g$-space coincide with variances
calculated in $\beta$-space.

Imposing similar restrictions based on higher order moments,
one might fully eliminate nonuniqueness.
%We are working on finding such restrictions.

Estimation of variances is another source of problems.
Formula (\ref{eq:Disper}) is of theoretical
importance, as it demonstrates that we have enough information
to estimate variances. However, it hardly can be used for numerical
computations as it involves {\em differences} of values
that we can only approximately estimate.
We are working on finding a better way to estimate variances.

%=======================================================================

\section{Conclusion}
\label{sec:Conclusion}

We developed a novel approach to analysis of categorical data
based on considering distribution laws of observed random variables
as realizations of another random variable.
This starting point leads to a fruitful development.

In the present article, we were able to obtain system of equations
(\ref{eq:MainEqSys}) and establish its properties in theorems
\ref{th:Main1}--\ref{th:Main4}.
This provides a base for an efficient numerical procedure
that gives (one form of) an answer to General GoM Problem.

We also believe that the approach in general,
and our results regarding system (\ref{eq:MainEqSys}) in particular,
may be successfully applied in other domains of statistics,
especially in latent structure analysis.

%=======================================================================

\bibliographystyle{alpha}

\bigskip
\bigskip
\bigskip

{
\parbox{3in}{\textsc{Mikhail Kovtun},\\
            \textsc{Igor Akushevich},\\
            \textsc{Kenneth G. Manton}
            \smallskip
            \\
            Center for Demographic Studies\\
            Duke University\\
            2117 Campus Drive\\
            Durham, NC 27708\\
            E-mail:  {\normalfont \url{mkovtun@cds.duke.edu}}\\
   \phantom{E-mail:} {\normalfont \url{aku@cds.duke.edu}}\\
   \phantom{E-mail:} {\normalfont \url{kmg@cds.duke.edu}}}
\hfill
\parbox{3in}{\textsc{H. Dennis Tolley}
            \medskip
            \\
            Department of Statistics\\
            Brigham Young University\\
            230 TMCB\\
            Provo, UT 84602\\
            E-mail: {\normalfont \url{tolley@byu.edu}}}
}

\end{document}